\documentclass[11pt,notitlepage, reqno]{amsart}

\usepackage{graphicx}
\usepackage{amsmath, amsfonts, amssymb, amsthm}
\usepackage[myheadings]{fullpage}
\usepackage[utf8]{inputenc}
\usepackage[margin=1in]{geometry}
\usepackage{enumitem}
\usepackage{xcolor}
\usepackage{caption}
\usepackage{tabularx}
\usepackage{fancyhdr}
\usepackage{hyperref}
\usepackage{tikz}
\usetikzlibrary{angles,quotes}
\usepackage{mdframed}
\usepackage{lipsum}
\usepackage{bm}
\usepackage{chngcntr}
\usepackage{imakeidx}
\usepackage{tkz-euclide}
\usepackage{blindtext}
\usepackage[inline]{asymptote}
\usepackage{mathdots}
\usepackage{tikz-cd,mathtools}
\usepackage{ mathrsfs }
\usepackage{float}
\restylefloat{table}

\usepackage{pgfplots, pgfplotstable}
\usepackage{multicol}
\usepackage{endnotes}
\usepackage{idxlayout}
\usepackage{xcolor, forest}
\usepackage{tocvsec2}
\usepackage[toc]{glossaries}
\usepackage{venndiagram}
\usepgfplotslibrary{statistics}
\usepackage[colorinlistoftodos]{todonotes} 
\usepackage[normalem]{ulem}
\usetikzlibrary{shapes.geometric}
\usepackage{mathtools}

\newcommand{\E}{\mathbb{E}}

\newcommand{\Z}{\mathbb{Z}}

\numberwithin{equation}{section}

\newtheorem{thm}{Theorem}[section]

\newtheorem{cor}[thm]{Corollary}

\newtheorem{definition}[thm]{Definition}

\newtheorem{lem}[thm]{Lemma}

\newtheorem{prop}[thm]{Proposition}
\newtheorem{rem}[thm]{Remark}
\theoremstyle{definition}
\newtheorem{convention}{Convention}

\newcommand{\block}[4]{\begin{pmatrix} {#1} & {#2} \\ {#3} & {#4}\end{pmatrix}}

\newcommand{\be}{\begin{equation}}
\newcommand{\ee}{\end{equation}}
\newcommand{\bea}{\begin{eqnarray}}
\newcommand{\eea}{\end{eqnarray}}

\DeclareMathOperator{\tr}{tr}
\DeclareMathOperator{\sw}{swirl}

\makeatletter
\newcommand{\addresseshere}{%
  \enddoc@text\let\enddoc@text\relax
}
\makeatother

\begin{document}

\title{Limiting Spectral Distributions of Families of Block Matrix Ensembles}

\author{Teresa Dunn}
\address{University of California-Davis}
\email{mtdunn@ucdavis.edu}

\author{Henry L. Fleischmann}
\address{University of Michigan}
\email{henryfl@umich.edu}

\author{Faye Jackson}
\address{University of Michigan}
\email{alephnil@umich.edu}

\author{Simran Khunger}
\address{Carnegie Mellon University}
\email{skhunger@andrew.cmu.edu}

\author{Steven J. Miller}
\address{Williams College}
\email{sjm1@williams.edu}

\author{Luke Reifenberg}
\address{University of Notre Dame}
\email{lreifenb@nd.edu}

\author{Alexander Shashkov}
\address{Williams College}
\email{aes7@williams.edu}

\author{Stephen Willis}
\address{Williams College}
\email{sdw2@williams.edu}

\thanks{This work was supported in part by NSF Grant DMS1947438 and Williams College. We thank our colleagues from the 2021 SMALL REU and Arup Bose for helpful comments.}



\date{\today}

\begin{abstract}
    We introduce a new matrix operation on a pair of matrices, $\sw(A,X),$ and discuss its implications on the limiting spectral distribution. In a special case, the resultant ensemble converges almost surely to the Rayleigh distribution. In proving this, we provide a novel combinatorial proof that the random matrix ensemble of circulant Hankel matrices converges almost surely to the Rayleigh distribution, using the method of moments. 
\end{abstract}

\maketitle

\markboth{Dunn, Fleischmann, Jackson, Khunger, Miller, Reifenberg, Shashkov, Willis}{Random Matrix Theory}


\section{Introduction} \label{sec:introduction section}
Random matrix theory was used by Eugene Wigner as a mechanism for modeling the limiting behavior of the energy distribution of heavy nuclei. The states of individual heavy nuclei are difficult to determine using the Schrödinger Equation, so instead one can examine the eigenvalues of random matrices and thereby obtain information about the statistical behavior of the system, as done in \cite{PhysicsRMT}. 

The techniques from nuclear physics were later abstracted to ensembles of random matrices. The motivation for choice of ensemble corresponded to the properties of physical systems. For example, this was the motivation for studying ensembles of real symmetric matrices, self-adjoint matrices, and Gaussian Orthogonal Ensembles. Given the importance of studying eigenvalues to both physics (as in \cite{wigner_1951, dyson_1962}) and to other fields of mathematics such as analytic number theory (as in \cite{katz_sarnak_1999, katz_sarnak_1999_2}), the eigenvalue distribution of the ensemble is the focus of study.  

In general, it is rare to find a named, closed form limiting distribution of the eigenvalue distributions for a given ensemble of random matrices. For example, in the ensemble of Toeplitz matrices studied in \cite{hammond_miller_2005} and \cite{bryc_dembo_jiang_2006}, the distribution seemed to be approaching a Gaussian distribution, but there were Diophantine obstructions with the index combinatorics of the random variable entries in the matrices. These obstructions prevented the distribution from being a Gaussian distribution, and a closed form is still not known. Following these difficulties, an attempt to overcome the obstructions and increase symmetry was done by adding palindromicity, this is sufficient to guarantee almost sure convergence to the Gaussian distribution \cite{massey_miller_sinsheimer}. Many other related ensembles have been thoroughly investigated, for example in \cite{blockMatRecentPaper, bose_mitra_2002, HighlyPalindromic, KologluKoppMiller2011, miller_swanson_tor_winsor_2015, KologluKoppMiller2011, beckwith_luo_miller_shen_triantafillou_2015}.

In this paper, we formulate a new matrix operation, ``swirl,'' based on the symmetry of the concentric even matrix ensemble. An example of a matrix in this ensemble is 


\[
    \begin{pmatrix}
      x_2 &  x_1 & x_0  & x_3  & x_3 &  x_0  & x_1 & x_2  \\
      x_1 &  x_0 & x_3  & x_2  & x_2 &  x_3  & x_0  & x_1 \\
      x_0 &  x_3 & x_2  & x_1  & x_1 &  x_2  &  x_3 & x_0   \\
      x_3 &  x_2 & x_1  & x_0  & x_0 &  x_1  &  x_2 & x_3   \\ 
      x_3 &  x_2 & x_1  & x_0  & x_0 &  x_1  &  x_2 & x_3  \\
      x_0 &  x_3 & x_2  & x_1  & x_1 &  x_2  &  x_3 & x_0 \\
       x_1   & x_0  & x_3  & x_2  & x_2 &  x_3  &  x_0 & x_1\\
      x_2    &  x_1      & x_0  & x_3  & x_3 &  x_0  &  x_1 & x_2
    \end{pmatrix}.
\]
Notably, the $x_i$ are variables drawn independently from a probability distribution with mean zero, variance one, and finite higher moments.

We chose this ensemble with the hope that by increasing symmetry we would be able to obtain a closed form for the spectral distribution of the matrices in the ensemble. 

It is advantageous to understand such matrices in block matrix form, as evidenced by \cite{blockMatRecentPaper}. In this vein, we split the matrices in the concentric even ensemble into blocks or quadrants and defined the swirl operation using two $N \times N$ input matrices, $A$ and $X$, to create the larger block matrix of size $2N \times 2N$ corresponding to the concentric even matrix, where $A$ is the upper right quadrant and $X$ is the exchange matrix. That is, 

    \begin{equation}
        \sw(A, X) \ = \ \block{AX}{A}{XAX}{XA}.
    \end{equation}

In concentric even matrices, $A$ is a circulant Toeplitz matrix and $AX$ is a circulant Hankel matrix. We reduce studying the circulant even ensemble to studying circulant Hankel matrices with several theorems about the behavior of $\tr(\sw(A,X))$ in Section \ref{sec:swirl matrices section}. Hankel matrices arise in a multitude of applications across fields of mathematics and physics: differential equations, functional analysis, statistics, probability theory, control theory, and more (see \cite{PhysRevE.103.042213, Peller2006HankelOA, iohvidov1982}, for example). Their symmetry also makes them a heavily studied family in random matrix theory, as in 
\cite{PhysRevE.103.042213,BOURGET2021103}. Circulant Hankel matrices also happen to be even centrosymmetric matrices, which have additional specialized applications in physics, for example in \cite{diele_sgura_2003}.

In Section \ref{sec:circulant Hankel section} we characterize summands in terms of the number of repeated entries and compute the moments via combinatorial degree of freedom arguments. By these methods, we obtain a novel combinatorial proof showing that the limiting spectral distribution of the random matrix ensemble of circulant Hankel matrices converges almost surely to the symmetrized Rayleigh distribution  (for an earlier proof relying on direct computation, see \cite{bryc_dembo_jiang_2006}). As we discuss in Appendix \ref{sec:Appendix A section}, our methods are generally applicable to many random matrix ensembles. In particular, we have the following theorem.
\begin{thm}
    Let $\mu_{A,N}(x)$ be the empirical spectral measure of the $N \times N$ circulant Hankel random matrix ensemble populated by entries from a sequence of random variables $A$ from a distribution $p$ with mean $0$, variance $1$, and finite higher moments. Then,
    \begin{equation}
       \lim_{N \to \infty}\mu_{A,N}(x) \to |x|e^{-x^2} 
    \end{equation}
almost surely.
\end{thm}

Notably, $|x|e^{-x^2}$ is the symmetrized Rayleigh distribution, with many known applications to physics (see \cite{siddiqui_1962}).

\begin{figure}[h!]
    \centering
    \includegraphics{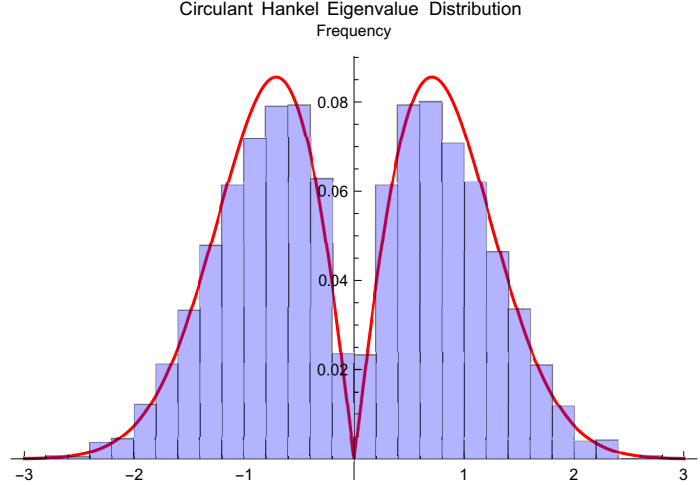}
    \caption{Histogram of eigenvalues for one hundred $40 \times 40$ random circulant Hankel matrices. A symmetrized Rayleigh distribution is shown in red.}
    \label{fig:my_label}
\end{figure}

Many block random matrix ensembles have been investigated in the past (for example, \cite{KologluKoppMiller2011}). Some of these have even yielded remarkably similar limiting empirical spectral distributions (see Figure 3 of \cite{miller_swanson_tor_winsor_2015}).


The swirl operation is very rich and lends itself to much further study. In particular, a natural next step is to study matrix ensembles determined by different choices of $A$ and $X$. We discuss some natural next steps in Section \ref{sec:future work section}.

%

\section{Preliminaries} \label{sec:preliminaries section}
We characterize the distribution of the eigenvalues of several random matrix ensembles by defining a spectral measure over subfamilies of random matrices from the ensemble. Let $A$ be an element of a family of $N \times N$ random matrices from some ensemble where the entries are drawn from a probability distribution $p$ with mean 0, variance 1, and finite higher moments. 

We use the Eigenvalue Trace Lemma to relate the eigenvalues to the matrix elements.


\begin{lem}[Eigenvalue Trace Lemma]\label{lem:evalue trace lemma}
    Let $\lambda_i(A)$ be the eigenvalues of an $N \times N$ matrix A. Then
    \begin{equation}
        \sum_{i=1}^N \lambda_i^k(A) \ = \ {\rm \tr}(A^k).
    \end{equation}
\end{lem}


Let $c$ be the number of eigenvalues of $A$, excluding those that are trivially zero. This is derived using rank arguments and is fixed for a given $N$ and matrix structure.

Then, we define the empirical spectral measure of $A$ as the following measure.
\begin{definition} \label{defn: spectral measure}
Let $p$ be a probability density function with mean $0$, variance $1$, and finite higher moments. Let $A$ be a family of $N \times N$ matrices with entries drawn independently from $p$. Then
    \begin{equation}
       \mu_{A,N}(x)dx \ := \ \frac{1}{c}\sum_{i = 1}^c\delta \left( x - \frac{\sqrt{c}\lambda_i(A)}{N} \right)dx,  
    \end{equation}
where $\delta(x)$ is the Dirac-delta functional, the $\lambda_i(A)$ are the nonzero eigenvalues of $A$, and $c$ is the number of eigenvalues in $A$ that are not trivially $0$.
\end{definition}
\begin{rem}
 The $\sqrt{c}/N$ scaling factor is derived heuristically from the Central Limit Theorem. By computing the trace of $A^2$ via the Eigenvalue Trace Lemma, we get
    \begin{equation}
     \mathbb{E}[\tr(A^2)] \ = \ N^2 \ = \ \sum_{i=1}^N \E[\lambda_i(A)^2],    
    \end{equation}
    suggesting that the magnitude of the eigenvalues must be roughly $N/\sqrt{c}$ each in expectation since the expectation of an entry squared is 1, by our definition of $p$.
\end{rem}

Via the method of moments, we will be able to understand the spectral distribution of these eigenvalues. In this instance, the convergence of the  moments of the spectral distribution is enough to show the convergence of the spectral distribution. From the definition of the spectral measure $\mu_{A,N}(x)$ in terms of the Dirac-delta functional, we may compute its moments.

\begin{rem}\label{rem: formula for moments}
    The moments of the spectral measure of $A$ are
    \begin{equation}
      M_k(A,N) \ := \ \int_{-\infty}^{\infty} x^k \mu_{A, N}(x) dx = \frac{c^{k/2 - 1}}{N^k} \sum_{i=1}^c \lambda_i^k(A).    
    \end{equation}
\end{rem}

Notice that by the Eigenvalue Trace Lemma, $M_k(A, N) = \frac{c^{k/2 - 1}}{N^k} {\rm \tr}(A^k).$

Finally, we're interested in averaging these moments over the entire family of matrices that $A$ belongs to. As is standard, we define the following.

\begin{definition}
    Let $M_k(N)$ be the average of $M_k(A,N)$ over all $A$ in our chosen family of matrices. 
\end{definition}

Our main result is that $\lim_{N \to \infty} M_k(N)$ exists  and that there is a universal limiting distribution for several families of matrices.

In the following work, we will calculate $M_k$ as a sum of terms (via the Eigenvalue Trace Lemma). We will show that some terms are negligible by showing that they are $O_k(1)$ (where $f(n) = O_k(g(n))$ if, for $k \in \mathbb{Z}^+$ fixed there exists $n_0, c$ such that for all $n > n_0, f(n) \leq g(n) + c)$.


In this paper, we investigate swirl ensembles and circulant Hankel matrices.

\begin{definition}
    An $N \times N$ circulant Hankel matrix $H_n = (a_{ij})$ is defined by the link relation $a_{ij} = a_{k\ell} \iff i+j \equiv k + \ell \pmod{N}$.
    \[ 
    J \ = \ \begin{pmatrix}
      b_0      & b_1      & b_2    & \cdots & b_{N-3}   & b_{N-2}   & b_{N-1}      \\
      b_1      & b_2      & b_3    & \cdots & b_{N-2}   & b_{N-1}       & b_0      \\
      b_2      & b_3      & b_4    & \cdots & b_{N-1}       & b_0         & b_2      \\
      \vdots   & \vdots   & \vdots &        & \vdots    & \vdots    & \vdots \\ 
      b_{N-3}  & b_{N-2}      & b_{N-1}      & \cdots & b_{N-6}     & b_{N-5}         & b_{N-4}      \\
      b_{N-2}      & b_{N-1}      & b_0      & \cdots & b_{N-5}     & b_{N-4}         & b_{N-3}   \\
      b_{N-1}      & b_0      & b_1      & \cdots & b_{N-4}         & b_{N-3}         & b_{N-2}   \\
    \end{pmatrix}, a_{ij}=b_{[i+j]_N}.
    \]
\end{definition}

Note that circulant Hankel matrices are the product of circulant Toeplitz matrices and exchange matrices, with the former considered in \cite{hammond_miller_2005, bryc_dembo_jiang_2006}.

\section{Swirl Matrices} \label{sec:swirl matrices section}
\subsection{Motivation}
The swirl operation was inspired by radially symmetric matrices of the following form:
\[
    \begin{pmatrix}
      x_2 &  x_1 & x_0  & x_3  & x_3 &  x_0  & x_1 & x_2  \\
      x_1 &  x_0 & x_3  & x_2  & x_2 &  x_3  & x_0  & x_1 \\
      x_0 &  x_3 & x_2  & x_1  & x_1 &  x_2  &  x_3 & x_0   \\
      x_3 &  x_2 & x_1  & x_0  & x_0 &  x_1  &  x_2 & x_3   \\ 
      x_3 &  x_2 & x_1  & x_0  & x_0 &  x_1  &  x_2 & x_3  \\
      x_0 &  x_3 & x_2  & x_1  & x_1 &  x_2  &  x_3 & x_0 \\
       x_1   & x_0  & x_3  & x_2  & x_2 &  x_3  &  x_0 & x_1\\
      x_2    &  x_1      & x_0  & x_3  & x_3 &  x_0  &  x_1 & x_2
    \end{pmatrix}.
\]
We refer to such matrices as ``concentric even matrices." Note that not only are the circles about the center of the matrix composed of equal entries, but also these entries are repeated in later circles such that each matrix entry appears an equal number of times. This was intentional in an effort to increase symmetry and derive a closed form limiting spectral distribution. For a $2N \times 2N$ matrix of this form, each entry appears exactly $4N$ times ($N$ times in each $N \times N$ quadrant). Upon close inspection, it is apparent that the $N \times N$ submatrix in the top right of a $2N \times 2N$ concentric even matrix is an $N \times N$ circulant Toeplitz matrix (which is not necessarily symmetric). Moreover, the other three quadrants of the matrix may be generated from this circulant Toeplitz matrix via a clockwise rotation of the entries. Indeed, for $A$ an $N \times N$ circulant Toeplitz matrix, and $J$ the $N \times N$ exchange matrix with 1's on the antidiagonal and zeroes elsewhere, the $2N \times 2N$ concentric even matrix is given by the following:
\[
    \block{AJ}{A}{JAJ}{JA}.
\]

This block decomposition of the concentric even matrices motivates the following definition and the focus of this section.

\begin{definition}  \label{defn: swirl}
    Let $A$ and $X$ be $N \times N$ matrices. We define {$\sw(A, X)$}
    as the $2N \times 2N$ matrix where 
    \begin{equation}
        \sw(A, X) \ = \ \block{AX}{A}{XAX}{XA}.
    \end{equation}
\end{definition}

We aim to characterize the limiting spectral distribution of $\sw(A,X)$. To do so, we relate $\tr(\sw(A,X)^k)$ to $\tr((AX)^k)$ via the Eigenvalue Trace Lemma.
\begin{rem}
    Observe that
    \begin{equation} \label{eqtn: expansion}
        \sw(A,X) \ = \  \block{AX}{A}{X A X}{X A} \ = \ \block{AX}{0}{0}{0} + \block{0}{A}{0}{0} + \block{0}{0}{X A X}{0} + \block{0}{0}{0}{X A}. 
    \end{equation}
\end{rem}
This observation vastly simplifies the computation of $\tr(((\sw(A,X))^k)$.

\begin{convention} \label {conv: block notation}
 We adopt a convenient shorthand notation for block matrices with four $N \times N$ blocks which are $0$ in 3 blocks. For example, a $2N \times 2N$ matrix of the form with zeroes necessarily everywhere except the top right corner will be referred to as a matrix $B_{12}.$ That is,
    $B_{12}$ is of the form
\[
    \block{0}{Y}{0}{0}
\]
for $Y$ an $N \times N$ matrix. Define $B_{11}$, $B_{21}, $ and $B_{22}$ similarly with the indices corresponding to the block that is not necessarily zero everywhere. 
\end{convention}
\begin{rem} \label{rem: block prod}
$B_{ij} B_{k\ell} = 0$ if $j \neq k$ and is of the form $B_{i\ell}$ otherwise.
\end{rem}

\subsection{Computing $\tr((\sw(A,X))^k)$}
Recall the following facts:
\begin{equation}
 \tr(CD) \ = \ \tr(DC)  
\end{equation}
and 
\begin{equation}
\tr(C +D) \ = \ \tr(C) + \tr(D) 
\end{equation}
for $N \times N$ matrices $C$ and $D.$ 
We are now ready to relate $\tr((\sw(A,X))^k)$ to $\tr((AX)^k)$.

\begin{thm} \label{thm: trSwirl thm}
For $A$ and $X$ both $N \times N$ matrices,  $\tr((\sw(A, X))^k) \ = \ 2^k\tr((AX)^k)$.
\end{thm}
\begin{proof}
 Moreover, any term in the expansion of 
 \begin{equation}
    (\sw(A,X))^k \,=\, \left(\block{AX}{0}{0}{0} + \block{0}{A}{0}{0} + \block{0}{0}{X A X}{0} + \block{0}{0}{0}{X A}\right)^k 
    \,=\, \sum B_{i_1 j_1}B_{i_2 j_2} \cdots B_{i_kj_k}
 \end{equation}
 is of the form $B_{i_1 j_k}$. Since trace is additive, we have that a term contributes 0 to the trace of $(\sw(A,X))^k$ if $i_1 \neq j_k$; if they are not equal, the main diagonal of the matrix is all zeroes.

As such, by Remark \ref{rem: block prod} and the above, the nonzero summands of 
$(\sw(A,X))^k$
correspond to products $B_{i_1 j_1}B_{i_2 j_2} \cdots B_{i_kj_k}$ where $j_\ell = i_{\ell+1}$ for $1 \leq \ell \leq k-1$ and $j_k = i_1$. There are $2^k$ such summands since one can choose the first indices of the $k$ matrices in the summand in $2^k$ ways. Then, the second indices are exactly determined by the above requirements.

Observe that the only nonzero block of a $B_{i_1j_k}$ matrix in the trace expansion of $\tr((\sw(A,X))^k)$ is a product of matrices. By the construction of swirl, this product begins with $A$ if $i=1$, and ending with $A$ if $j=2$.  Also observe that this product begins with an $X$ if $i=2$, and ends with $X$ if $j=1$. All such products will start with one of $A$ or $X$ and end with the opposite. These products will also not have consecutive repeated $A$'s or $X$'s. These properties follow from Remark \ref{rem: block prod} and the definition of swirl. 

In order for this product of matrices to contribute to the trace, note that the first and last index of such a product must be equal (or else it will not be a diagonal entry). Thus, there must be an equal number of matrices of the form $B_{12}$ and $B_{21}$ in any contributing product. As such, each nonzero summand in the expansion of $\tr((\sw(A, X))^k)$ (with $\sw(A, X)$ expressed as $B_{11} + B_{12} + B_{21} + B_{22}$) will be of the form  $\tr((XA)^k)$ or $\tr((AX)^k)$. Consequently, there are $2^k$ such nonzero contributing terms and
\begin{equation}
    \tr((\sw(A, X))^k) \ = \ 2^k\tr((AX)^k),
\end{equation}
 by the cyclic property of trace.
\end{proof}


Given that the trace of the $k$\textsuperscript{th} power of a matrix completely determines the $k$\textsuperscript{th} moment of its empirical spectral distribution, Theorem \ref{thm: trSwirl thm} allows us to reduce characterizing the limiting spectral distribution of $\sw(A,X)$ ensembles to characterizing the limiting spectral distribution of $AX$ matrices.

\subsection{Iterating $\sw$}
Another interesting avenue for swirl is iterating the operation.

\begin{definition} \label{def: it swirl}
Let $A,X$ be $N \times N$ matrices. Let $X_{k}$ be the block matrix with $2^{k-1}$ $X$'s on the anti-diagonal and zeroes elsewhere. Note that $X_1=X.$ Then, set
\begin{equation}
    \sw^k(A,X) \ := \ \sw(\ldots \sw(\sw( \sw(A,X_1), X_2), X_3), \ldots), X_k)
\end{equation}
where $\sw$ is repeated $k$ times in the above.
\end{definition}
We begin by analyzing the trace of iterated swirl matrices.
\begin{prop} \label{prop: more swirls}
Fix $A,X$ both $N \times N$ matrices such that $X^2=I$ and $k$ a nonnegative integer. Then
    \begin{equation}
        \tr(\sw^k(A,X)) \ = \ 2^k \tr(AX).
    \end{equation}
\end{prop}
\begin{proof}
We prove by induction. For $k = 1$, this follows from Theorem \ref{thm: trSwirl thm}. Now, assume this holds for $r-1$ for $r \geq 2$. Then, 
\begin{align}
     \sw^r(A,X) \ = \ \sw(\sw^{r-1}(A,X), X_{r}) \\
    \ = \ \block{X_{r}\sw^{r-1}(A,X)}{\sw^{r-1}(A,X)}{X_{r} \sw^{r-1}(A,X) X_{r}}{\sw^{r-1} (A,X) X_{r}} \nonumber.
\end{align}
This implies
\begin{align}
    \tr(\sw^r(A,X)) \ = \ \tr(X_{r} \sw^{r-1}(A,X)) + \tr( \sw^{r-1}(A,X) X_{r}) \\
    \ = \ 2\tr( X_{r} \sw^{r-1}(A,X)) \nonumber.
\end{align}
Now
\begin{align}
   X_{r} \sw^{r-1}(A,X) \ = \ \block{0}{X_{r-1}}{X_{r-1}}{0} \block{X_{r-1}\sw^{r-2}(A,X)}{\sw^{r-2}(A,X)}{X_{r-1}\sw^{r-2}(A,X)X_{r-1}}{\sw^{r-2}(A,X)X_{r-1}} \\
   \ = \ \block{\sw^{r-2}(A,X)X_{r-1}}{X_{r-1}\sw^{r-2}(A,X)X_{r-1}}{\sw^{r-2}(A,X)}{X_{r-1}\sw^{r-2}(A,X)} \nonumber.
\end{align}

Thus 
\begin{align}
    \tr(X_{r}\sw^{r-1}(A,X)) &= \tr(\sw(\sw^{r-2}(A,X), X_{r-1})) \\
    \ = \ \tr(\sw^{r-1}(A,X)) \nonumber \\
    \ = \ 2^{r-1}\tr(AX) \nonumber
\end{align}
by induction.

Therefore,
\begin{equation}
    \tr(\sw^r(A,X)) \ = \ 2^r\tr(AX).
\end{equation}
The result then follows by induction.
\end{proof}

\begin{rem} \label{rmk: block repetition iterated swirl}
Alternatively, observe that, since $X^2 = I$, $\sw^{\ell}(A,X)$ is just the block matrix of $\sw(A,X)$ repeated $4^{\ell-1}$ times. This means $\tr(\sw^{\ell}(A,X)) = 2^{\ell-1}\tr(\sw(A,X)) = 2^{\ell}\tr(AX)$.
\end{rem}

If we wish to study the moments of ensembles of such matrices, we need to understand the trace of powers of the iterated swirl matrices. We reduce this to an analysis of $\tr((AX)^k)$ in the following proposition.
\begin{prop} \label{prop: powers iterated swirl}
Fix $A,X$ to be $N \times N$ matrices such that $X^2=I$ and $k$ and  $l$  nonnegative integers.  Then
    \begin{equation}
        \tr((\sw^{\ell}(A,X))^k) \ = \ 2^{k\ell} \tr((AX)^k).
    \end{equation}
\end{prop}
\begin{proof}
We prove by induction on $\ell$. For $\ell = 1$, this follows from Theorem \ref{thm: trSwirl thm}. 
Now, assume 
\begin{equation}
    \tr((\sw^{\ell}(A,X))^k) \ = \ 2^{k\ell} \tr((AX)^k)
\end{equation}
holds for $\ell = r \geq 1$. We show it holds for $\ell = r+1$.
By Definition \ref{def: it swirl}, 
\begin{equation}
    \sw^{\ell+1}(A,X) \ = \ \sw(\sw^{\ell}(A, X_{\ell}), X_{\ell+1}).
\end{equation}
So, 
\begin{align*}
    \tr((\sw^{\ell+1}(A,X))^k) \ = \ \tr((\sw(\sw^{\ell}(A, X), X_{\ell+1}))^k) \\
    \ = \ 2^k \tr((\sw^{\ell}(A, X)X_{\ell+1})^k)\nonumber \,
\end{align*}
with the last step following from Theorem \ref{thm: trSwirl thm}.

Let $B = \sw^{\ell-1}(A, X)$. Then, 
\begin{align}
    \sw^{\ell}(A, X)X_{\ell+1} \ = \ \block{BX_\ell}{B}{X_\ell BX_\ell}{X_\ell B}  \block{0}{X_\ell}{X_\ell}{0} \\ 
    \ = \ \block{BX_\ell}{B}{X_\ell BX_\ell}{X_\ell B} = \sw^{\ell}(A, X) \nonumber
\end{align}
with the last step following from the assumption that $X^2 = I$.

Therefore 
\begin{align*}
    \tr((\sw^{\ell+1}(A,X))^k) &= 2^k \tr((\sw^{\ell}(A, X)X_{\ell+1})^k) \\
    &= 2^k \tr((\sw^{\ell}(A, X))^k) \nonumber\\
    &= 2^{k(\ell+1)}\tr((AX)^k)\nonumber \,
\end{align*}
with the last step from the inductive hypothesis. The result then follows by induction.
\end{proof}

\subsection{The Product of Swirl and its Transpose}

If we assume that $X$ is a permutation matrix, then $\tr(\sw(A,X)\sw(A,X)^T)$ reduces to understanding $\tr(AA^T)$. This is a useful quantity to understand if $A$ and $X$ are chosen such that $\sw(A,X)$ does not necessarily have real eigenvalues.

\begin{prop} \label{prop: transpose swirl}
    Fix $A,X$ to be $N \times N$ matrices with $X$ a permutation matrix.  Then
    \begin{equation}
        \tr((\sw(A,X)\sw(A,X)^T)^k) \ = \ 2^{2k}\tr((AA^T)^k).
    \end{equation}
\end{prop}

\begin{proof}
    Let $S$ = $\sw(A,X)\sw(A,X)^T$.  We show by induction that
    
    \begin{equation}
        S^k\ =\ 2^{2k-1}\block{(AA^T)^k}{(AA^T)^kX^T}{X(AA^T)^k}{X(AA^T)^kX^T}.
    \end{equation}
    
    For the base case consider $S^k$ for $k=1$.  We have
    \begin{equation}
        S^1\ =\ \block{AX}{A}{X A X}{X A} \block{(AX)^T}{(XAX)^T}{A^T}{(XA)^T}.
    \end{equation}
    
    This yields 
    
    \begin{equation}
        =\block{(AX)(AX)^T+AA^T}{(AX)(XAX)^T+A(XA)^T}{(XAX)(AX)^T+(XA)A^T}{(XAX)(XAX)^T+(XA)(XA)^T}.
    \end{equation}
    
    Expanding the transpose terms yields 
    
    \begin{equation}
        =\block{AXX^TA^T+AA^T}{AXX^TA^TX^T+AA^TX^T}{XAXX^TA^T+XAA^T}{XAXX^TA^TX^T+XAA^TX^T}.
    \end{equation}
    
    Recall $X$ is a permutation matrix $XX^T=I$. Thus, we have 
    
    \begin{equation}
        =2\block{AA^T}{AA^TX^T}{XAA^T}{XAA^TX^T}.
    \end{equation}
    
    Now assume that the inductive hypothesis holds for $k=n$; we will show it holds for $k=n+1$.  Rewrite $S^{n+1}$ as $S S^n$. Then
    
    \begin{equation}
        S^{n+1} = 2\block{AA^T}{AA^TX^T}{XAA^T}{XAA^TX^T} 2^{2n-1}\block{(AA^T)^n}{(AA^T)^nX^T}{X(AA^T)^n}{X(AA^T)^nX^T},
    \end{equation}
    
    by induction. Matrix multiplication yields
    
    \begin{equation}
        =2^{2n} \block{(AA^T)^{n+1}+(AA^T)^{n+1}}{(AA^T)^{n+1}X^T+(AA^T)^{n+1}+X^T}{X(AA^T)^{n+1}+X(AA^T)^{n+1}}{X(AA^T)^{n+1}X^T + X(AA^T)^{n+1}X^T}.
    \end{equation}
    
    Simplifying we have
    
    \begin{equation}
        =2^{2n+1}\block{(AA^T)^{k+1}}{(AA^T)^{k+1}X^T}{X(AA^T)^{k+1}}{X(AA^T)^{k+1}X^T}.
    \end{equation}
    
    This completes the inductive argument.
    
    Now calculating the trace is trivial.  Note by the cyclic property of trace, $\tr(S^k) = 2^{2k-1} \tr((AA^T)^k) + 2^{2k-1} \tr(X(AA^T)^kX^T)) = 2^{2k} \tr((AA^T)^k)$.  
\end{proof}
Here the limiting spectral distribution reduces to the a scaled semi-circle distribution, which is handled in \cite{wigner_1958}.

\subsection{Limiting Spectral Distribution of Swirled Matrix Ensembles} \label{subsec: reduce to lim spec dist.}
From the previous work in this section, we can reduce the analysis of swirl matrix ensembles to the analysis of matrix product ensembles. We consider the empirical spectral measure defined in Definition \ref{defn: spectral measure}. In this case, from Remark \ref{rmk: block repetition iterated swirl}, for $A$ and $X$ both $N \times N$ matrices, and $\ell \geq 1$, $\sw^\ell(A,X)$ has the same number of trivial nonzero eigenvalues, $c$, as $\sw^\ell(A,X)$. Let $B_{N2^\ell} = \sw^\ell(A, X)$. Then the empirical spectral measure of $B_{N2^\ell}$ is given by 
\begin{equation}
    \mu_{A,N2^\ell}(x)dx \ := \ \frac{1}{c}\sum_{i = 1}^{c}\delta \left( x - \frac{\sqrt{c}\lambda_i(A)}{N2^l} \right)dx. 
\end{equation}
From Definition \ref{defn: spectral measure} and Proposition \ref{prop: powers iterated swirl}, the $k$\textsuperscript{th} moment of the spectral distribution in this case is thus
\begin{equation}
    \frac{2^{\ell k}c^{k/2 -1 }}{2^{\ell k}N^{k+1}} \E[\tr((AX)^k)] \ = \ \frac{c^{k/2 -1}}{N^{k+1}}\E[\tr((AX)^k)],
\end{equation}
which does not depend on $\ell$. As such, the limiting spectral distribution of swirl is the same for any number of iterations, $\ell$.

\section{Circulant Hankel Matrices} \label{sec:circulant Hankel section}

In all the ensembles that follow, we assume that the matrices are constructed from a sequence of independently and identically distributed random variables (i.i.d.r.v.) with distribution $p$ having mean 0, variance 1, and finite higher moments. We assign elements of this sequence to matrix entries according to the symmetry of our given ensemble. 

\subsection{Moments via powers of $AX$}
From Theorem \ref{thm: trSwirl thm}, studying the trace of the even concentric swirl matrices reduces to studying the trace of powers of $H_N = A_NJ_N$, with $H_N$ the $N \times N$ circulant Hankel matrix, $A_N$ the $N \times N$ circulant Toeplitz matrix, and $J_N$ the $N \times N$ exchange matrix. The matrix ensemble of circulant Hankel matrices is exceptional in its own right; its limiting spectral distribution converges almost surely to a symmetrized Rayleigh distribution (as shown in \cite{bryc_dembo_jiang_2006}). In this section, we provide a new combinatorial proof of this remarkable result. We begin by defining the \textit{empirical spectral measure} for this ensemble of matrices. This measure, for the normalized eigenvalues of our matrix $H$, is given by the following definition.

\begin{definition} \label{defn: spectral measure Hankel}
The empirical spectral measure of a random $N \times N$ circulant Hankel matrix is
\begin{equation}
\mu_{H_N}(x)\text{dx} \ := \ \frac{1}{N}\sum_{i = 1}^N\delta \left( x - \frac{\lambda_i(H_N)}{\sqrt{N}} \right)\text{dx}.    
\end{equation}
where $\delta(x)$ is the Dirac-delta functional and the $\lambda_i$ are the non-zero eigenvalues of $H_N$. 
\end{definition}
\begin{rem}
The $\sqrt{N}$ scaling factor is derived heuristically. By computing the trace of $H_N^2$, we obtain
\begin{equation}
\mathbb{E}[\tr(H_N^2)] \,=\, N^2 \,=\, 
\sum_{i=1}^N \lambda_i(H_N)^2,    
\end{equation}
suggesting that the eigenvalues must be roughly $\sqrt{N}$ each in expectation.
\end{rem}

In order to use the method of moments, we compute the $k$th moment for the empirical spectral distribution of a random matrix $H_N$, $\mu_{H_N}(x)$.

\begin{rem} \label{rem: kth mom}
The $k$th moment of the empirical spectral distribution of the random matrix $H_N$, averaged over an ensemble, is given by
\begin{equation}
M_k(N) \ := \ \int_{- \infty}^{\infty}x^k \mu_{H_N}(x) dx \ = \ \frac{1}{N^{\frac{k}{2}+ 1}}\sum_{i=1}^N\E[\lambda_i^k(H_N)] \ = \ \frac{1}{N^{\frac{k}{2}+ 1}} \E[\tr(H_N^k)].    
\end{equation}
We use $M_k$ to denote $\lim_{N \to \infty} M_k(N)$.
\end{rem}

This standard computation follows from the properties of the Dirac delta functional and the Eigenvalue Trace Lemma.

\begin{prop} \label{ prop: moments 1 + 2 G_N}
We have 
$M_1 = 0$ and $M_2 = 1$.
\end{prop}
\begin{proof}
The first moment is immediate from $\mathbb{E}[\tr(H_N)] = 0$. The second moment follows from substituting $\mathbb{E}[\tr(H_N^2)] = N^2$ into the formula in Remark \ref{rem: kth mom}.
\end{proof}

In order to compute $M_k$ for $k \geq 2$ we consider the limiting behavior of the terms in the sum combinatorially. It is useful to note the following fact.

\begin{rem} \label{ rem: indexing HankCirc}
In $H_N = A_N J_N$, $c_{ij} = c_{k\ell}$ if and only if $i + j \equiv_N k + \ell$, where we index the matrix beginning at 0.
\end{rem}

We begin by showing the odd moments of the limiting spectral distribution are all zero.

\begin{thm} \label{thm: oddMomZero}
We have $M_{2k+1} \ = \ 0$ for all $k \in \mathbb{Z}_{\geq 0}$.
\end{thm}
\begin{proof}
First, we analyze $\tr(H_N^{2k+1})$ via the Eigenvalue Trace Lemma. Observe that since the entries of $H_N$ are independent, if any are to the first power in a summand in the expansion of $\tr(H_N)$, the expected value of the entire summand is zero. For example, 
\begin{equation}
\E[h_{i_1i_2}h_{i_2i_3}h_{i_3i_1}] \,=\, \E[x_ax_bx_b] \,=\, \E[x_a]\E[x_b^2] \,=\, 0 \cdot 1 \,=\, 0.        
\end{equation}

 
 Thus, at a minimum, the entries must be matched in pairs with at least one triple in all of the contributing terms. We bound the number of ways to construct such summands. There are at most $k$ distinct entries $x_{j_i}$ in a given contributing summand by this pairing argument. We can choose such entries in less than $N^k$ ways. Then, we can specify the matrix index of one of the terms in the summand in $N$ ways. There are at most $k^k$ ways to assign each factor in the summand to a particular $x_{j_i}$ and then the choice of one index of a matrix entry completely determines the remaining matrix indices via Remark \ref{ rem: indexing HankCirc}. This implies there are at most $k+1$ degrees of freedom for any choice of grouping (and the number of ways to assign factors in the summand to $x_{j_i}$ is $O_k(1)$). So, the number of contributing summands is $O_k(N^{k+1})$. Note that each grouping of $n_i$ matrix entries equal to $x_{j_i}$ contributes 
\begin{equation}
    \mathbb{E}[(x_{j_i})^{n_i}] \,=\, p_{n_i} \,=\,  O_k(1)   
\end{equation}
since $p$ has finite higher moments by assumption. As such, each contributing term contributes $O_k(1)$ to $M_k$.

Substituting into our formula from Remark \ref{rem: kth mom}, we get:
\begin{equation}
 M_{2k+1}(N) \,=\, \frac{O_k(N^{k+1})}{N^{\frac{2k+1}{2} + 1}} \,=\, O_k(N^{-1/2}).    
\end{equation}
Taking the limit as $N \to \infty$, we achieve the desired result.
\end{proof}

Next, we show $M_{2k} = k!$ for all $k$ and thus $\lim_{N \to \infty} \mu_{H_N}(x)$, averaged over all $H_N$ converges to the symmetrized Rayleigh distribution.

We begin with a sample calculation showing $M_4 = 2$ to build intuition for the proof.

\begin{prop} \label{prop: G4thMoment}
We have $M_4 = 2$.
\end{prop}
\begin{proof}
Note that
\begin{equation}
\tr(H_N^4) \ = \ \sum_{i_1=1}^N\sum_{i_2=1}^N\sum_{i_3=1}^N\sum_{i_4=1}^N h_{i_1i_2} h_{i_2i_3} h_{i_3i_4} h_{i_4i_1}   
\end{equation}
where $h_{i_ji_{j+1}}$ is the matrix entry of $H_N$ at the $i_j$\textsuperscript{th} row and and $i_{j+1}$\textsuperscript{th} column.
As before, if any of the random variables in a summand is not equal to any of the others, we can write the expectation of the whole summand as a product of the expectation of the singleton term and the rest of the summand by the independence of our random variables. Since all of the random variables have mean 0, such a term contributes zero. As such, there are only two options for contributing summands: four equal matrix entries or two pairs of equal matrix entries.

\begin{description}
\item[\textit{Case 1}] In this case, there are four summands that are all matched. That is, $c_{i_1i_2} = c_{i_2i_3} = c_{i_3i_4} = c_{i_4i_1}.$
Up to relabeling, the first case yields the system of equations
\begin{align}
    i_1 \,+\, i_2 \,&\equiv_N\, i_2 \,+\, i_3 \\
    i_2 \,+\, i_3 \,&\equiv_N\, i_3 \,+\, i_4 \nonumber \\ 
    i_3 \,+\, i_4  \,&\equiv_N\, i_4 \,+\, i_1 \nonumber \\
    i_4 \,+\, i_1 \,&\equiv_N\, i_1 \,+\, i_2 \nonumber
\end{align}
This implies $i_1 \equiv_N i_3$ and $i_2 \equiv_N i_4$, leaving only two free variables. Since there are only 2 degrees of freedom in this case and each of our i.i.r.d. random variables have finite moments by assumption, terms of this kind contributes $O(N^2)$ to the expectation of $\tr(H_N^4)$. Thus, by Remark \ref{rem: kth mom}, such terms contribute $\lim_{N \to \infty}\frac{O(N^2)}{N^3} = 0$ to the fourth moment in the limit.

Notably, the system of equations corresponds to the equation matrix
\[
\begin{pmatrix}
1 & 1 & -1 & -1 \\
0 & 1 & 0 & -1 \\
-1 & 0 & 1 & 0 \\
0 & -1 & 0 & 1 
\end{pmatrix}
\]
which has nullity 2. Thus, since vectors satisfying this system of equations are exactly those in the null space of this matrix, there are $O(N)$ valid linear combinations of basis vectors of the null space, and the random variables have finite fourth moments, such terms contribute $O(N^2)$ to the expectation of $\tr(H_N^4)$. This alternate linear algebraic formulation is used in our the proof of Theorem \ref{thm: even moments G_N}.

\item [\textit{Case 2}] In this case, all summands are paired.
This case of matching the random variables into pairs has two subcases.
\begin{description}
\item[\textit{Subcase 2.1}]
Pair nonadjacent random variables, that is, $c_{i_1i_2} = c_{i_3i_4}, c_{i_2i_3} = c_{i_4i_1}$. This pairing yields the following system of equations:
\begin{align}
    i_1 \,+\, i_2 \,&\equiv_N\, i_3 \,+\, i_4 \\
    i_2 \,+\, i_3 \,&\equiv_N\, i_4 \,+\, i_1 \nonumber\\
     i_3 \,+\, i_4 \,&\equiv_N\, i_1 \,+\, i_2 \nonumber\\
     i_4 \,+\, i_1 \,&\equiv_N\, i_2 \,+\, i_3 \nonumber
\end{align}
This implies $i_2 \equiv_N i_4$ and $i_1 \equiv_N i_3$. Thus there are only two degrees of freedom in this case and it does not contribute in the limit.

Note that the equation matrix corresponding to the system of equations has nullity 2, an alternative proof that this case cannot contribute.
\item[\textit{Subcase 2.2}]
Pair adjacent random variables. For example, $c_{i_1i_2} = c_{i_2i_3}$ and $c_{i_3i_4} = c_{i_4i_1}$. Note that there are two such pairings. Up to relabeling, this pairing yields the following system of equations:
\begin{align}
    i_1 \,+\, i_2 \,&\equiv_N\, i_2 \,+\, i_3 \\
    i_2 \,+\, i_3 \,&\equiv_N\, i_1 \,+\, i_2 \nonumber\\
    i_3 \,+\, i_4 \,&\equiv_N\, i_4 \,+\, i_1 \nonumber\\
    i_4 \,+\, i_1 \,&\equiv_N\, i_3 \,+\, i_4 \nonumber
\end{align}
This implies only $i_1 \equiv_N i_3$, yielding 3 degrees of freedom. Thus, the terms in this case contribute in the limit. Fixing $i_1$, there is a unique choice for $i_3$ and $N$ choices for both $i_2$ and $i_4$, yielding $2N^3$ choices total after iterating over all $i_1$ and both choices of pairing orientation.

Substituting into Remark \ref{rem: kth mom}, we then get that this case contributes precisely $2$ to $M_4$ in the limit and 
$M_4 = 2$,
since this is the only contributing case.
\end{description}
\end{description}
\end{proof}

We see that only a select few of the summands in the computation of even moments contribute in the limit. We formalize this observation in the following lemmas.

\begin{lem} \label{lem: evenMomentsPaired}
For even moments $M_{2k}$, where $k \geq 1$, the only contributing summands $x_{j_1}^{n_1} \cdots x_{j_\ell}^{n_\ell}$ in the trace expansion  are those where $n_i = 2$ for all $1 \leq i \leq \ell$.
\end{lem}
\begin{proof}
Consider any summand in $\tr(H_N^{2k})$, $x_{j_1}^{n_1} \cdots x_{j_\ell}^{n_\ell}$, where 
\begin{equation}
  \sum_{i=1}^{\ell}n_i \,=\, 2k  
\end{equation}
and each $n_i \geq 1$.
Now, if any $n_i = 1$, the expectation of the summand is 0. So, we may assume each $n_i \geq 2$. If there is at least one factor in the summand with $n_r \geq 3$, there are less than $k+1$ degrees of freedom of terms with such groupings—there are at most $k-1$ ways to choose the $x_{j_i}$ and an additional $N$ ways to fix a matrix index in some term (which then induces a constant in $N$ number of possible arrangements).

Each such term contributes 
\begin{equation}
\lim_{N \to \infty} \frac{O(N^k)}{N^{k+1}} \,=\, 0    
\end{equation}
to $M_{2k}$, as desired.
\end{proof}
\begin{rem} \label{rem: degFreedomMax}
The argument in Theorem \ref{thm: oddMomZero} and Lemma \ref{lem: evenMomentsPaired} also shows that there are at most $k+1$ degrees of freedom when assigning $x_{j_i}$ in pairs in the computation of $\tr(H_N^{2k})$. 
\end{rem}

For the following arguments, consider an index ``even" if its subscript is even. Similarly define odd indices. For example, consider
\begin{equation}
\tr(H_N^4) \,=\, \sum_{i_1=1}^N\sum_{i_2=1}^N\sum_{i_3=1}^N\sum_{i_4=1}^N h_{i_1i_2}h_{i_2i_3}h_{i_3i_4}h_{i_4i_1}.    
\end{equation}
We view $i_1, i_3$ as the ``odd" indices and $i_2,i_4$ as the ``even" indices. 
\begin{lem} \label{lem: odd-even pairings}
The pairings of odd indices to even indices contribute $k!$ to $M_{2k}$.
\end{lem}
\begin{proof}
Consider the system of equations resultant in this case. Each relation can be assumed to be of the general form
\begin{equation}
i_j 
\,+\, i_{j+1} \,\equiv_N\, i_\ell \,+\, i_{\ell+1}    
\end{equation}

for $j$ even and $\ell$ odd. Note in particular that all even indices arise on the left hand side of such relations as the first term and all odds similarly as the first term on the right hand side. Since in such relations each index is added to the subsequent index, every index appears in a sum exactly once on both sides of the equations.

Interpret these equations as $1 \times N$ row vectors with ones in the entries corresponding to the indices on the left hand side of the relations and negative ones to those on the right hand, as in Proposition \ref{prop: G4thMoment}. Now, from the above observation, the sum of these $k$ row vectors is 0. This implies they are linearly dependent. This means the matrix given by this system of equations has nullity at least $k+1$. Note that vectors $\bm{x}$ in the null space are exactly solutions to 
\begin{equation}
E\bm{x} \,=\, 0    
\end{equation}
for $E$ the matrix of these row vectors. This implies we have $k+1$ degrees of freedom in this case. From Remark $\ref{rem: degFreedomMax}$, we thus have exactly $k+1$ degrees of freedom, so these pairings contribute exactly their constant term to $M_{2k}$ in the limit.

To count the number of odd-even pairings, we choose an odd and an even index to pair in $k^2$ ways. Then we repeat until there are no indices left to pair, yielding $(k!)^2$. However, we introduced an arbitrary ordering on the pairs in this process, so we correct by dividing by $k!$, yielding $k!$ as desired. Note that, given the choice of a single index and pairings, every index is determined uniquely (regardless of the modulo $N$).
\end{proof}

Now we complete the proof by showing that the other pairings of indices do not contribute in the limit.

\begin{thm} \label{thm: even moments G_N}
$M_{2k} \ = \ k!$. 
\end{thm}
\begin{proof}
From Lemma \ref{lem: odd-even pairings}, it suffices to show that any arrangement of pairs including an odd-odd or even-even index matching will not contribute. One way to do so is to show that the $k$ row vectors corresponding to the resultant system of equations are all linearly independent and thus the rank of the corresponding matrix is $k$, implying a nullity of $k$ and less than $k+1$ degrees of freedom. Note that there being an odd-odd index pairing implies that there must be an even-even index pairing. 
\begin{description}
\item[Step 1] We will show that if there is an odd-odd index pairing then the equations corresponding to even-odd index pairings are linearly independent as row vectors.

Fix a relation given by such an even-odd index pairing. Each side of each relation of the form $i_\ell + i_{\ell+1} \equiv_N i_r + i_{r+1}$ can be conceptualized as as a ``first" index (matched index) plus a ``second" index. For the sake of consistency, when converting such relations into row vectors (by moving all the terms to a single side), we negate the side with the odd first index. In order to show linear independence of the even-odd row vectors, it suffices to show that no nonempty linear combination of them sums to 0.

Note that each index appears at most twice amongst the odd-even pair relations. In particular, if we fix 
\begin{equation}
 i_\ell \,+\, i_{\ell+1} \,\equiv_N\, i_{\ell + 2i + 1} \,+\, i_{\ell + 2i + 2}   
\end{equation}
to be in our linear combination, with $\ell$ even, this yields a row vector of the form
\[
(0, 0, \ldots, 1, 1, 0, \ldots, 0, -1 -1, 0, \ldots, 0)
\]
with the $1$'s in the $\ell$ and $(\ell+1)$\textsuperscript{st} positions and the $-1$'s in the $(\ell+2i+1)$\textsuperscript{st} and $(\ell+2i+2)$\textsuperscript{th}  positions and all other positions $0$ (note the indices range from $1$ to $2k$).

Crucially, in order to form a linear combination of even-odd vectors summing to 0, we must nullify each index in the sum. Each index appears at most once as a first and a second index. Given our signing convention, even indices are positive as first indices and negative as second indices. Odds indices are negative as first indices and positive as second indices. Note that, since each index occurs at most once as a first index and a second index, the two expressions cannot be exactly equal. In order to cancel out the positive contribution of $i_\ell$ to the $\ell$\textsuperscript{th} column, we need to add the term including $i_{\ell}$ as a second index. However, we then must cancel out the contribution of $i_{\ell-1}$ as a first index by including it as a second index. To do that, we must include $i_{\ell-2}$ as a first index. As such we see that in order to cancel out the contribution of each necessary term, we need to include every term as both a first and second index. However, by assumption, there is an odd-odd index. Therefore, not every index has a row vector corresponding to it as a first and second index. Thus, we cannot cancel out the contribution to every column and there is no nonempty linear combination of vectors corresponding to the even-odd pair expressions that equals 0. We conclude that these row vectors are linearly independent.

\item[Step 2] We will show that the row vectors corresponding to odd-odd pairs of indices cannot be part of any nonempty linear combination of row vectors summing to zero (the proof follows for even-even pairs as well). 

Suppose indices $i_r$ and $i_{r+ 2i}$ are paired for $r$ odd and $i \geq 1$.
The corresponding row vector is of the form 
\[
(0, 0, \ldots, 1, 1, 0, \ldots, 0, -1 -1, 0, \ldots, 0)
\]
with $1$'s in the $r$\textsuperscript{th} and $(r+1)$\textsuperscript{st} indices and $-1$'s in the $(r+2i)$\textsuperscript{th} and $(r+2i+1)$\textsuperscript{st} indices. As before, to cancel out the contribution of $i_r$, we need $i_r$ to appear as a second index and contribute negatively. As a result, $i_{r-1}$ must appear as a first index and contribute negatively. Then we need $i_{r-1}$ to appear as a second index and contribute positively to cancel out that contribution. This requires $i_{r-2}$ to appear as a first index and contribute positively. However, this implies that all first odd indices must contribute positively and all first even indices must contribute negatively to achieve total cancellation. We know this cannot be the case as there is an odd-odd pair and one of the first odd indices must thus contribute negatively. As such, odd-odd and even-even pairs cannot be a part of linear combinations of the row vectors summing to 0.

We conclude that a linearly dependent family of row vectors must be a subset of the even-odd pair row vectors if it exists. However, from Step 1, this is impossible. So, all of the row vectors are linearly independent. Thus, if there are odd-odd or even-even index pairs, the rank of the matrix is $k$ and the nullity is $k$. Since the nullity of this matrix is a upper bound on the degrees of freedom in this case, such pairings will not contribute in the limit. The only remaining pairings are all odd-even. From Lemma \ref{lem: odd-even pairings} we conclude $M_{2k} = k!$.
\end{description}
\end{proof}

We use in Theorem \ref{thm: rayleigh} that the moments of the limiting spectral distribution of the Circulant Hankel ensemble are the same as the moments of a symmetrized Rayleigh distribution. More broadly, a Rayleigh distribution is a Weibull distribution with fixed parameters. For our purposes, denote the Weibull distribution with scale parameter $\lambda$ and shape parameter $k$ by the following: 
\begin{equation}
 f(x; \lambda, k) \,=\, \frac{kx^{k-1}}{x\lambda^k}e^{-(x/\lambda)^k},   
\end{equation}
for $x \geq 0$ and $0$ otherwise. As our eigenvalue distributions are symmetric, we symmetrize the distribution by replacing $x$ with $|x|$ and dividing through by 2 to retain 
\begin{equation}
 \int_{-\infty}^{\infty} \frac{f(|x|; \lambda, k)}{2}dx \,=\, 1.   
\end{equation}

This symmetrization notably has no effect on the even moments of the distribution and zeroes all the odd moments.

The $(2n)$\textsuperscript{th} moment of a Weibull distribution $f(x; \lambda k)$ is given by
\begin{equation}
m_{2n} \,=\, \lambda^{2n}\Gamma(2n/k + 1).    
\end{equation}

When this distribution is Rayleigh, i.e., $k = 2$ and $\lambda = 1$, the $n$\textsuperscript{th} moment for $n$ even is then $n!$.

\begin{thm} \label{thm: rayleigh}
Let $A$ denote an infinite sequence of values drawn from a distribution $p$ with mean $0$, variance $1$, and finite higher moments. As $N \to \infty$ the limiting spectral measure of the circulant Hankel random matrix ensemble converges almost surely to the limiting spectral distribution given by the $M_m$'s, the symmetrized Rayleigh distribution:
\begin{equation}
f(x) \,=\, |x|e^{-x^2}.    
\end{equation}
\end{thm}
\begin{proof}
This follows by the exact same argumentation as in Section 6 of \cite{hammond_miller_2005} with plus rather than minus modulo $N$ .
\end{proof}

We may then conclude that the limiting spectral distribution of iterated swirl ensembles on $A$ circulant Toeplitz and $J$ an exchange matrix also converges almost surely to a symmetrized Rayleigh distribution.

\begin{cor} \label{cor: it swirl dist}
 Let $G_{2N} = \sw(A,J)$ for $J$ the $N \times N$ exchange matrix and $A$ a random $N \times N$ circulant Toeplitz matrix. As $N \to \infty$, the limiting spectral measure of this ensemble converges almost surely to a symmetrized Rayleigh distribution.
\end{cor}
\begin{proof}
From  the observation that $\sw(A,J)$ trivially has half of its rows repeated, $G_{2N}$ has only $N$ nontrivial, nonzero eigenvalues.  The empirical spectral measure of the $2N \times 2N$ matrix $B_{2N}$  is thus given by the following equation:
\begin{equation}
\mu_{B_{2N}}(x)dx \,\coloneqq\, \frac{1}{N}\sum_{i = 1}^N\delta \left( x - \frac{\lambda_i(B_{2N})}{2\sqrt{N}} \right)dx.    
\end{equation}
See Definition \ref{defn: spectral measure} for the derivation of the scaling factor. From Theorem \ref{thm: trSwirl thm}, the $k$\textsuperscript{th} moment of the limiting spectral distribution of this ensemble equals 
\begin{equation}
\lim_{N \to \infty} \frac{1}{N^{k/2 +1}}\E[\tr(H_N^k)].    
\end{equation}
As such, the $k$\textsuperscript{th} moment in this case is exactly the $k$\textsuperscript{th} moment of the limiting spectral distribution of $H_N$. The result then follows from Theorem \ref{thm: rayleigh}.
\end{proof}

\begin{rem}
Note that, for an ensemble such that $\sw(A,X)$ has no repeated rows, $\sw(A,X)$ would not have the same limiting spectral distribution as $AX$. Indeed, its moments would be $2^k$ times the moments of the limiting spectral distribution of $AX$. In the case of a Weibull distribution, this would only increase the $\lambda$ scaling parameter of the Weibull distribution by a factor of $\sqrt{2}$.
\end{rem}

\section{Future Work} \label{sec:future work section}

The obvious next step is to study broader matrix ensembles related to $\sw(A,X)$. A good starting point is ensembles with $X^2 = I$, due to the following theorem of Tao and Yasuda \cite{tao_david_2002_1236140}.
 
\begin{thm} [Tao-Yasuda \cite{tao_david_2002_1236140}, 2002]
Let $A$ and $X$ be real symmetric matrices with $X^2 = I$. 
\begin{itemize}
    \item $AX \ = \ XA$ if and only if the spectrum of $A$ equals the spectrum of $XA$ up to sign.
    \item $AX \ = \ -XA$ if and only if the spectrum of $A$ equals the spectrum of $XA$ multiplied by $i$.
\end{itemize}
\end{thm}

In particular, if we choose ensembles $A$ and $X$ such that $A$ and $X$ are $N \times N$ real symmetric matrices, $X^2 = I$ and $AX = XA$, then $AX$  has all real eigenvalues. 

Another interesting direction is study the even powers of non-symmetric swirl ensembles. Proposition \ref{prop: transpose swirl} provides a useful starting point for such investigations.

Finally, given that circulant Toeplitz and circulant Hankel matrices yield rare named, closed form limiting spectral distributions, it seems likely that they possess some intrinsic, special properties. Inspired by the work of \cite{bose_basak_2018, banerjee_bose_2011, bose_hazra_saha_2011, miller_swanson_tor_winsor_2015} on matrices with patterns governed by link functions, we investigated circulant matrices with link functions along different diagonals, but found the results disappointingly uninteresting. These results are summarized in Appendix \ref{sec:Appendix A section}.

\appendix 

\section{} \label{sec:Appendix A section}

When computing the moments of the limiting empirical spectral measures of our ensembles we converted our problem of finding degrees of freedom of contributing summands in the trace to a problem of calculating the nullity of a matrix. As a specific example, we can calculate the $n$th moment of the Hankel ensemble by looking at set partitions of $\{1,2, \ldots n\}$ and calculating the rank of  matrices of the form 
\begin{equation}
M_\pi\,=\,(I_n - P_\pi)B_n(1,1),   
\end{equation}
where $I_n$ is the identity matrix, $P_\pi$ is a permutation matrix, and $B_n(s,t)$ is a matrix with $s$ in the diagonal and $t$ to the right of the diagonal. The matrix $B_4(s,t)$ is written below as an example.
\begin{equation}
B_4(s,t) \,=\, \begin{bmatrix}
s & t & 0 & 0 \\
0 & s & t & 0 \\
0 & 0 & s & t \\
t & 0 & 0 & s \\
\end{bmatrix}    
\end{equation}
The permutation matrix $P_{\pi}$ corresponds to the particular matchings of indices in the summand corresponding to equal matrix entries. The nullity of the matrix $M_{\pi}$ gives the degrees of freedom of assignments of entries to groups that contribute in the case $\pi$. Iterating over all set partitions and substituting into the formula in Definition \ref{rem: formula for moments}, we can easily show certain configurations do not contribute in the limit.

We generalize by considering matrices which are constant along certain circulant lines/diagonals. We call these $(s,t)$-ensembles. Formally, an $N \times N$ matrix is in the $(s,t)$-ensemble if 
\begin{equation}
si + tj \,\equiv_N\, sk + t\ell \,\implies\, a_{ij} \,=\, a_{k\ell}.    
\end{equation}

An example of a matrix within the (1,2)-ensemble is:
\begin{equation}
A_4(1,2)\,=\,\begin{bmatrix}
x_1 & x_2 & x_3 & x_4\\
x_3 & x_4 & x_1 & x_2\\
x_1 & x_2 & x_3 & x_4\\
x_3 & x_4 & x_1 & x_2
\end{bmatrix}    
\end{equation}
Note how the equivalent entries, denoted here by entries of the same value, are spaced apart by $s=1$ row movements and $t=2$ column movements. This creates the appearance of a matrix where every $s$ rows, the entries are horizontally permuted by $t$ columns. From this, we obtain the idea of slope as we describe the relationship between equivalent entries in $(s,t)$-ensemble patterns. This generalization is the idea of a polynomial link function in the literature, except now modulo $N$ (see \cite{bose_basak_2018, banerjee_bose_2011, bose_hazra_saha_2011, miller_swanson_tor_winsor_2015}).

However, this generalization is insufficient if we intend to study symmetric matrices. Thus, we strengthen our condition to 
\begin{equation}
 a_{ij} \ = \ a_{mn} \,\iff\, si+tj \,\equiv_N\, sm+tn \ \text{ OR }\ ti \,+\, sj \,\equiv_N\, tm \,+\, sn.   
\end{equation}
This allows us to generalize the special behavior of both circulant Hankel and circulant Toeplitz matrices. 

Below we is a $4 \times 4$ matrix that  belonging to the $(1,1)$-ensemble
\[\begin{bmatrix}
x_1 & x_2 & x_3 & x_4\\
x_2 & x_3 & x_4 & x_1\\
x_3 & x_4 & x_1 & x_2\\
x_4 & x_1 & x_2 & x_3
\end{bmatrix}.\]
These represent the circulant Hankel matrices. Moreover,  $4 \times 4$ matrices from the $(1,-1)$-ensemble are of the form
\[\begin{bmatrix}
x_1 & x_2 & x_3 & x_4\\
x_4 & x_1 & x_2 & x_3\\
x_3 & x_4 & x_1 & x_2\\
x_2 & x_3 & x_4 & x_1
\end{bmatrix}.\]
These represent the circulant Toeplitz matrices.

Inspired by the fact that both the circulant Hankel and circulant Toeplitz matrices admitted, we generalize the structure of these matrices in the hopes of finding a broader class of matrices with limiting empirical spectral distributions given by named probability distributions.  To this end, notice how the elements of the aforementioned families cascade through the matrix with "slope" $\pm 1$.  It is this notion of slope which we wish to generalize, and will be made more concrete in what follows.

With this change in parameters, the same pattern of equivalent entries being $s$ rows and $t$ columns away persists, and the main observable change is in the number of equivalence classes of matrix entries that appear.

The elements at the indices generated by $(s,t)$ and $(t,s)$ form a group: \begin{equation}
 H \,=\, \left< (s,t), (t,s) \right> \leq (\Z/N\Z) \times (\Z/N\Z).    
\end{equation}

Notice that $H$ is a normal subgroup of $(\Z/N\Z) \times (\Z/N\Z)$. We want to understand the number of cosets associated to $((\Z/N\Z) \times (\Z/N\Z))/H,$ as this will give us the number of unique elements within an $N\times N$ matrix with the aforementioned rule.

In the simplest cases we get the matrices we studied in the main portion of the paper.  When $s$ and $t$ both equal $1$, the resultant ensemble is Hankel, with the number of cosets increasing as $N$ increases. It is this positive slope that reflects the symmetry of the matrix that is lacking in the ciruclant Toeplitz. Likewise, if $s$ and $t$ are units with opposite signs, i.e., $s=1$ and $t=-1$, the resulting matrices are all Toeplitz.  Similarly, the number of cosets increases consistently with $N$.  The number of cosets is important to consider because it indicates the amount of variation within the matrix, the more cosets there are, the fewer zero eigenvalues appear.

Now, as we vary $s$ and $t$, new patterns arise in the family of matrices and consequently the number of cosets. This variance is a function of the positioning of equivalent entries. With these new $s$ and $t$ values, the spacing between the placement of the entries changes, and there are some very interesting patterns to the numbers of cosets and the qualities of symmetry. However, among all these patterns, it appears that the only ones that remain symmetric are circulant Hankel. Besides those, we continue to observe circulant Toeplitz matrices appearing at certain intervals and numbers of cosets.

For $N$ coprime to $s$ and $t$, we observe that when 
\begin{equation}
s \,\equiv\, t \pmod{[H\,:\,(\Z/N\Z) \times (\Z/N\Z)]},    
\end{equation}
where $[H\,:\,(\Z/N\Z) \times (\Z/N\Z)]$ is the index of $H$ in $(\Z/N\Z) \times (\Z/N\Z)$, the matrices yielded are circulant Toeplitz. Alternatively, when 
\begin{equation}
s \,\equiv\, -t \pmod{[H\,:\,(\Z/N\Z) \times (\Z/N\Z)]}    
\end{equation}
the matrices yielded are circulant Hankel.

However, we found that whenever we consider $N \times N$ matrices from the $(s,t)$-ensemble  with $s, t \neq \pm 1$, then the limiting spectral distribution is uninteresting.  This is because the number of cosets for a matrix in this ensemble appears bounded by a constant times $\text{gcd}((s+t)(s-t), N)$. Indeed, the matrix becomes a block matrix with many repetitions of a much smaller matrix, deferring its spectral distribution to that smaller matrix ensemble. When $s$ and $t$ are units (up to sign), we find computationally that the number of cosets is proportional to $\text{gcd}(0, N)$, which is just $N$.  Because the number of cosets is proportional to $N$, the number of eigenvalues grows as we increase the size of the matrix.  However, in the other case, number of nonzero eigenvalues is fixed, preventing a new distribution from arising.


\bibliographystyle{alpha}

\bibliography{references}

\newcommand{\etalchar}[1]{$^{#1}$}
\begin{thebibliography}{MSTW15}

\bibitem[BB11]{banerjee_bose_2011}
Sayan Banerjee and Arup Bose.
\newblock Noncrossing partitions, catalan words, and the semicircle law.
\newblock {\em Journal of Theoretical Probability}, 26(2):386–409, 2011.

\bibitem[BB18]{bose_basak_2018}
Arup Bose and Anirban Basak.
\newblock Balanced toeplitz and hankel matrices.
\newblock {\em Patterned Random Matrices}, page 131–142, 2018.

\bibitem[BBV{\etalchar{+}}19]{blockMatRecentPaper}
Keller Blackwell, Neelima Borade, Charles VI, Noah Luntzlara, Renyuan Ma,
  Steven Miller, Mengxi Wang, and Wanqiao Xu.
\newblock Distribution of eigenvalues of random real symmetric block matrices.
\newblock 2019.

\bibitem[BDJ06]{bryc_dembo_jiang_2006}
Włodzimierz Bryc, Amir Dembo, and Tiefeng Jiang.
\newblock Spectral measure of large random hankel, markov and toeplitz
  matrices.
\newblock {\em The Annals of Probability}, 34(1), 2006.

\bibitem[BG21]{PhysRevE.103.042213}
Eugene Bogomolny and Olivier Giraud.
\newblock Statistical properties of structured random matrices.
\newblock {\em Phys. Rev. E}, 103:042213, 2021.

\bibitem[BHS11]{bose_hazra_saha_2011}
Arup Bose, Rajat~Subhra Hazra, and Koushik Saha.
\newblock Convergence of joint moments for independent random patterned
  matrices.
\newblock {\em The Annals of Probability}, 39(4), 2011.

\bibitem[BLM{\etalchar{+}}15]{beckwith_luo_miller_shen_triantafillou_2015}
Olivia Beckwith, Victor Luo, Steven~J Miller, Karen Shen, and Nicholas
  Triantafillou.
\newblock Distribution of eigenvalues of weighted structured matrix ensembles.
\newblock {\em Integers}, 15:A21, 2015.

\bibitem[BM02]{bose_mitra_2002}
Arup Bose and Joydip Mitra.
\newblock Limiting spectral distribution of a special circulant.
\newblock {\em Statistics \& Probability Letters}, 60(1):111–120, 2002.

\bibitem[Bou21]{BOURGET2021103}
A.~Bourget.
\newblock Spectral distribution of families of hankel matrices.
\newblock {\em Linear Algebra and its Applications}, 624:103--120, 2021.

\bibitem[DS03]{diele_sgura_2003}
F.~Diele and I.~Sgura.
\newblock Centrosymmetric isospectral flows and some inverse eigenvalue
  problems.
\newblock {\em Linear Algebra and its Applications}, 366:199--214, 2003.

\bibitem[Dys62]{dyson_1962}
Freeman~J. Dyson.
\newblock The threefold way. algebraic structure of symmetry groups and
  ensembles in quantum mechanics.
\newblock {\em Journal of Mathematical Physics}, 3(6):1199–1215, 1962.

\bibitem[FM09]{PhysicsRMT}
Frank Firk and Steven Miller.
\newblock Nuclei, primes and the random matrix connection.
\newblock {\em Symmetry}, 1, 2009.

\bibitem[HM05]{hammond_miller_2005}
Christopher Hammond and Steven~J. Miller.
\newblock Distribution of eigenvalues for the ensemble of real symmetric
  toeplitz matrices.
\newblock {\em Journal of Theoretical Probability}, 18(3):537–566, 2005.

\bibitem[JMP10]{HighlyPalindromic}
Steven Jackson, Steven Miller, and Thuy Pham.
\newblock Distribution of eigenvalues of highly palindromic toeplitz matrices.
\newblock {\em Journal of Theoretical Probability}, 25, 2010.

\bibitem[KKM11]{KologluKoppMiller2011}
Murat Koloğlu, Gene~S. Kopp, and Steven~J. Miller.
\newblock The limiting spectral measure for ensembles of symmetric block
  circulant matrices.
\newblock {\em Journal of Theoretical Probability}, 26(4):1020--1060, 2011.

\bibitem[KS99a]{katz_sarnak_1999}
Nicholas Katz and Peter Sarnak.
\newblock Random matrices, frobenius eigenvalues, and monodromy.
\newblock {\em Colloquium Publications}, 1999.

\bibitem[KS99b]{katz_sarnak_1999_2}
Nicholas~M. Katz and Peter Sarnak.
\newblock Zeroes of zeta functions and symmetry.
\newblock {\em Bulletin of the American Mathematical Society}, 36(01):1–27,
  1999.

\bibitem[MMS07]{massey_miller_sinsheimer}
Ada Massey, Steven~J. Miller, and John Sinsheimer.
\newblock Distribution of eigenvalues of real symmetric palindromic toeplitz
  matrices and circulant matrices.
\newblock {\em Journal of Theoretical Probability}, 20(3):637--662, 2007.

\bibitem[MSTW15]{miller_swanson_tor_winsor_2015}
Steven~J. Miller, Kirk Swanson, Kimsy Tor, and Karl Winsor.
\newblock Limiting spectral measures for random matrix ensembles with a
  polynomial link function.
\newblock {\em Random Matrices: Theory and Applications}, 04(02), 2015.

\bibitem[Pel06]{Peller2006HankelOA}
V.~Peller.
\newblock Hankel operators and their applications.
\newblock {\em IEEE Transactions on Automatic Control}, 51:383--385, 2006.

\bibitem[SGT82]{iohvidov1982}
Iohvidov~Iosif Semënovič., Israel Gohberg, and Gerard Philip~Antoine.
  Thijsse.
\newblock {\em Hankel and Toeplitz matrices and forms: algebraic theory}.
\newblock Birkhäuser, 1982.

\bibitem[Sid62]{siddiqui_1962}
M~M Siddiqui.
\newblock Some problems connected with rayleigh distributions.
\newblock {\em Journal of Research of the National Bureau of Standards, Section
  D: Radio Propagation}, 66D(2):167, 1962.

\bibitem[TY02]{tao_david_2002_1236140}
David Tao and Mark Yasuda.
\newblock {A Spectral Characterization of Generalized Real Symmetric
  Centrosymmetric and Generalized Real Symmetric Skew-Centrosymmetric
  Matrices}.
\newblock {\em SIAM Journal on Matrix Analysis and Applications}, 2002.

\bibitem[Wig51]{wigner_1951}
Eugene~P. Wigner.
\newblock On the statistical distribution of the widths and spacings of nuclear
  resonance levels.
\newblock {\em Mathematical Proceedings of the Cambridge Philosophical
  Society}, 47(4):790–798, 1951.

\bibitem[Wig58]{wigner_1958}
Eugene~P. Wigner.
\newblock On the distribution of the roots of certain symmetric matrices.
\newblock {\em The Annals of Mathematics}, 67(2):325, 1958.

\end{thebibliography}

\addresseshere

\end{document}